\documentclass[11pt, reqno]{amsart}
\usepackage{times,amsmath,amssymb,amsthm,amscd, color}
\usepackage[all]{xy}
\usepackage{xcolor}

\newtheorem{thm}{Theorem}[section]
\newtheorem{cor}[thm]{Corollary}
\newtheorem{lem}[thm]{Lemma}
\newtheorem{pro}[thm]{Proposition}

\theoremstyle{definition}

\newtheorem{rmk}[thm]{Remark}

\numberwithin{equation}{section}

\newcommand{\R}{{\mathbf R}}
\newcommand{\C}{{\mathbf C}}
\newcommand{\Cn}{{\mathbf C^n}}
\newcommand{\Bn}{{\mathbf B_n}}
\newcommand{\D}{{\mathbf D}}
\newcommand{\re}{{\text{Re}}\,}

\DeclareMathOperator*{\supp}{Supp}

\title
[Estimates of the exponential Bergman kernel on the unit ball]
{Pointwise estimates of the Bergman kernel with an exponential weight on the unit ball}

\author{Hong Rae Cho, Soohyun Park*}


\begin{document}

\begin{abstract}
We consider the weighted Bergman space $A^2_\psi(\Bn)$
of all holomorphic functions on $\Bn$ square integrable with respect to 
a particular exponential weight measure $e^{-{\psi}} dV$ on $\Bn$, where 
\begin{align*}
\psi(z):=\frac{1}{1-|z|^2}.
\end{align*}
We prove the following estimate for the Bergman kernel $K_\psi(z,w)$ of $A^2_\psi(\Bn)$:
\begin{align*}
	|K_\psi(z,w)|^2\le C\frac{e^{\psi(z)+\psi(w)}}{{\rm Vol}(B_\psi(z,1)){\rm Vol}(B_\psi(w, 1))}e^{-\varepsilon d_\psi(z,w)},
\quad z, w\in\Bn,
\end{align*}
where $d_\psi$ is the Riemannian distance induced by the potential function $\psi$ and $B_\psi(z,1)$ is the $d_\psi$-ball
of center $z$ and radius $1$.
The result is motivated by Christ \cite{Chr}.
\end{abstract}

\maketitle

\section[Introduction]{Introduction}

Let $\C$ denote the complex field and $\Cn$ denote the cartesian product of $n$ copies of $\C$ for positive integer $n$. 
For $z=(z_1, \cdots , z_n)$ and $w=(w_1, \cdots , w_n) \in \Cn$, we define the Hermitian inner product
\begin{align*}
\langle z,w\rangle= \sum_{j=1}^n z_j \bar w_j,
\end{align*}
and the associated norm 
$|z| = \langle z,z\rangle^{1/2}. $
The open unit ball of $\Cn$ is denoted by $\Bn = \left\{z \in \Cn : |z| < 1\right\}.$
Let $dV$ be the Lebesgue volume measure on $\Cn$.
For a strictly plurisubharmonic function
\begin{align*}
\psi(z):=\frac{1}{1-|z|^2},
\end{align*}
let the weighted Bergman space with respect to $\psi$ 
\begin{align*}
  A^2_\psi(\Bn) := 
  \left\{ 
  f : \Bn \rightarrow \C \text{~holomorphic function} ; \int_\Bn \left| f(z) \right|^2 e^{-\psi(z)} dV(z)  < \infty 
  \right\}.
\end{align*}
It is a Hilbert space with the Hermitian inner product
\begin{align*}
\langle f, g\rangle_{\psi} := \int_{\Bn}f(z)\overline{g(z)} e^{-\psi(z)}\, dV(z)
\end{align*}
for $f, g \in A^2_{\psi}(\Bn)$.
Lemma \ref{SMVP} guarantees that each point evaluation $L_zf=f(z)$ is bounded on $A^2_\psi(\Bn)$.  
By the Riesz representation theorem, there is a holomorphic function $K_z \in A^2_{\psi}(\Bn)$ satisfying 
\begin{align*}
f(z) = \langle f, K_z\rangle_{\psi} = \int_{\Bn}f(w)\overline{K_z(w)} e^{-\psi(w)}\, dV(w).
\end{align*}
The function $K_{\psi}(z,w):=\overline{K_z(w)}$ is called the Bergman kernel for $A^2_{\psi}(\Bn)$ and
has the property $K_z(w)=\overline{K_w(z)}$. 

\bigskip

Function spaces with generalized weights have attracted a lot of attention in recent years (one can refer for Bergman type spaces: \cite{AH}, \cite{CP}, \cite{Dos}, \cite{GP}, 
\cite{HLS}, \cite{KM}, \cite{LR}, \cite{PP}, \cite{SV2}; for Fock type spaces: \cite{Chr}, \cite{Dal}, \cite{Del}, \cite{LR}, \cite{Lin}, \cite{MO}, \cite{SV1}, \cite{SY}) 
since new techniques,
different from those used for standard Bergman (or, Fock) spaces, are required. 

In 1991, Christ \cite{Chr} obtained pointwise upper bound of the Bergman kernel of $A_\varphi^2(\C)$ by means of the distance with $\Delta \phi$ for a large class of subharmonic function $\varphi$ under certain conditions, which was weakened by \cite{MO}. 
In the case of $\D$, similar estimates have been obtained in \cite{AH} and \cite{HLS}. For the case of $\Cn$, see \cite{Dal}, \cite{Del}, \cite{Lin}, \cite{SV1}, and \cite{SY}. 
For the case of $\Bn$, \cite{CP} and \cite{SV2} have presented estimates on weighted Bergman kernels.

The aim of this paper is to illustrate an upper bound of the Bergman kernel for $A^2_{\psi}(\Bn)$ in terms of a distance function with $\psi$. While \cite{CP} provide a diagonal estimate, our work shows an off-diagonal estimate, which is not covered by the results in \cite{SV2}. 

Our main result is the following.
	\begin{thm}
	There are constants $C>0$ and $0<\varepsilon<\sqrt{2}$ satisfying
	\begin{align*}
	|K_{\psi}(z,w)|^2 \leq C \frac{e^{\psi(z)}e^{\psi(w)}}{(1-|z|^2)^{2n+1} (1-|w|^2)^{2n+1}} e^{-\varepsilon d_\psi(z, w)} 
	\quad \text{for} \quad z, w \in \Bn.
	\end{align*}
	\end{thm}
Here, $d_\psi(z,w)$ is the Riemannian distance between $z$ and $w$ induced by the potential function $\psi$.
We note that ${\rm Vol}(B_\psi(z,r))\simeq (1-|z|^2)^{2n+1}$, where $B_\psi(z,r)$ is the $d_\psi$-ball
(see Corollary \ref{cor_vol}).

\bigskip

Throughout this paper, the notation $U(z)\lesssim V(z)$ (or equivalently $V(z)\gtrsim U(z)$) 
means that there is a constant $C$ such that $U(z) \le CV(z)$ holds for all $z$ in the set in question, which may be a space of functions or a set of numbers. If both $U(z)\lesssim V(z)$ and $V(z)\lesssim U(z)$, then we write $U(z)\simeq V(z)$.

\section{Hermitian geometry induced by the potential function $\psi$}

\subsection{Hermitian metric induced by $\psi$}
We recall that 
\begin{align*}
\psi(z)=\frac{1}{1-\lvert z \rvert^2},\quad z\in\Bn.
\end{align*}
The complex Hessian of $\psi$ is defined by 
\begin{align*}
H_{\psi}:=\left(\frac{\partial^2{\psi}}{\partial z_j\partial\bar z_k}\right)_{n \times n}.
\end{align*}

\begin{lem}[\cite{Har}, \cite {SM}]  \label{SM}
Let $S$ be an invertible $n\times n$ matrix, and $u, v\in\Cn$. Then
\begin{align*}
\det(S+uv^T)=\det(S)(1+v^TS^{-1}u)
\end{align*}
and, if $1+v^T S^{-1}u\neq 0$,
\begin{align*}
(S+uv^T)^{-1}=S^{-1}-\frac{S^{-1}uv^TS^{-1}}{1+v^TS^{-1}u}.
\end{align*}
\end{lem}

We introduce an auxiliary matrix 
\begin{align*}
A(z)=\left(z_j \bar z_k\right)_{n\times n}
\end{align*}
for each $z=(z_1, \cdots, z_n)\in\Cn$.
Let $P_z$ denote the orthogonal projection in $\Cn$ onto the complex line $\{\lambda z : \lambda\in\C\}$, where $z$ is an
arbitrary point in $\Cn\setminus \{0\}$. It will be convenient to let $P_0$ denote the identity transformation.
If we identify linear transformations on $\Cn$ with $n\times n$ matrices via the standard basis of $\Cn$ (so that 
the adjoint of a linear transformation is just the conjugate transpose of the corresponding matrix), then for $z\neq 0$,
\begin{align*}
A(z)=|z|^2P_z.
\end{align*}

\begin{lem}\label{lem for H} 
The complex Hessian $H_{\psi}$ has the following properties:
\begin{enumerate}
\item [(a)] 
$ H_{\psi}(z) = \frac{1}{(1-\lvert z \rvert^2)^3}\left((1-\lvert z \rvert^2)I_{n\times n} + 2\overline{A(z)}\right)$.
\item [(b)] 
$H_{\psi}(z)^{-1} = {(1-\lvert z \rvert^2)^2}\left(I_{n\times n} - \frac{2}{1+\lvert z \rvert^2}\overline{A(z)}\right)$.
\item [(c)] 
$\det H_{\psi}(z) = \frac{1+\lvert z \rvert^2}{\left(1-\lvert z \rvert^2\right)^{2n+1}}$.
\item [(d)] 
$H_{\psi}(z)= \frac{1+\lvert z \rvert^2}{(1-\lvert z \rvert^2)^3}\bar P_z + \frac{1}{(1-\lvert z \rvert^2)^2}\bar Q_z,$
where $Q_z = I-P_z$.
\end{enumerate}
\end{lem}
\begin{proof}
It follows that
\begin{align*}
\frac{\partial^2}{\partial z_j\bar\partial z_k}\left(\frac{1}{1-|z|^2}\right)
=\frac{\delta_{j,k}}{(1-|z|^2)^2}+\frac{2 \bar z_j z_k}{(1-|z|^2)^3}.
\end{align*}
This shows that
\begin{align*}
H_{\psi}(z) = \frac{1}{(1-\lvert z \rvert^2)^3}\left((1-\lvert z \rvert^2)I_{n\times n} + 2\overline{A(z)}\right).
\end{align*}
The facts $A(z)=\lvert z \rvert^2 P_z$ and $P_z + Q_z = I$ are used for the case of (d).
  
Let $u=(\bar z_1, \cdots, \bar z_n)$, $v= (z_1, \cdots ,z_n)$ be column vectors, then $u v^T=\left(\bar z_j z_k\right)_{jk}=\overline{A(z)}$ and $v^T u = \lvert z \rvert^2$.
From (a), it is obtained
\begin{align}\label{Hessian matrix reform}
H_{\psi}(z) = \frac{2}{(1-\lvert z \rvert^2)^3}\left(S(z) + u v^T \right)
\end{align}
where $S(z)=\frac{(1-\lvert z \rvert^2)}{2}I_{n\times n}$.
The Sherman-Morrison formula (Lemma \ref{SM}) gives
\begin{align*}
\left(S(z) + u v^T \right)^{-1} 
&= S(z)^{-1}-\frac{S(z)^{-1} uv^{T} S(z)^{-1}}{1+v^{T} S(z)^{-1} u}\\
&= \frac{2}{1-\lvert z \rvert^2} I_{n\times n} - \frac{4}{(1+\lvert z \rvert^2)(1-\lvert z \rvert^2)}\overline{A(z)}
\end{align*}
which shows (b).

Also, the Sherman-Morrison formula gives
\begin{align*}
\det \left(S(z) + u v^T \right) 
&= \left( 1 + v^{T} S(z)^{-1} u\right) \det S(z) \\
&= \frac{(1 + \lvert z \rvert^2)(1-\lvert z \rvert^2)^{n-1}}{2^n}
\end{align*}
which provides (c) from \eqref{Hessian matrix reform}.
\end{proof}

\begin{rmk}
The statement (d) in Lemma \ref{lem for H} means that for $n \geq 2$ and $z \neq 0$, the matrix $ H_{\psi}(z)$ has $n$ eigenvalues, namely, an eigenvalue $\frac{1+\lvert z \rvert^2}{(1-\lvert z \rvert^2)^3}$ with eigenspace $[z]$, and $(n-1)$ eigenvalues $\frac{1}{(1-\lvert z \rvert^2)^2}$ with eigenspace ${\C^n} \ominus [z]$. 
\end{rmk}

The complex Hessian $H_\psi$ is positive definite and 
the metric
\begin{align}\label{hermitian}
h_{\psi}=\sum_{j,k=1}^n\left(\frac{\delta_{j,k}}{(1-|z|^2)^2}+\frac{2 \bar z_j z_k}{(1-|z|^2)^3}\right)dz_j\otimes d\bar z_k
\end{align}
defines a Hermitian metric on $\Bn$. 
The fundamental $(1,1)$-form
\begin{align*}
i\partial\bar\partial\psi(z)
=i\sum_{j,k=1}^n\left(\frac{\delta_{j,k}}{(1-|z|^2)^2}+\frac{2 \bar z_j z_k}{(1-|z|^2)^3}\right) dz_j\wedge d\bar z_k
\end{align*}
associated to $h_\psi$  defines a K\"ahler metric on $\Bn$.
In fact, $(\Bn, i\partial\bar\partial\psi)$ is a complete K\"ahler manifold. 
Let $\xi\in\Cn$. From (\ref{hermitian}) we have
\begin{align*}
|\xi|_{h_\psi}^2
&=h_\psi(\xi,\xi)\\
&=\frac{2|\langle\xi,z\rangle|^2}{(1-|z|^2)^3}+\frac{|\xi|^2}{(1-|z|^2)^2}\\
&=\frac{1+|z|^2}{(1-|z|^2)^3}|P_z\xi|^2+\frac{|Q_z\xi|^2}{(1-|z|^2)^2}.
\end{align*}
Let $(h_\psi)_{jk}$ be the components of the matrix $H_\psi$ and 
$h_\psi^{jk}$ the components of the inverse matrix $H_\psi^{-1}$. 
Then
\begin{align}\label{inverse}
h_\psi^{jk}=(1-|z|^2)^2\delta_{j,k}
-\frac{2(1-|z|^2)^2}{1+|z|^2}\bar z_j z_k.
\end{align}
Let $\alpha=\sum_{j=1}^n\alpha_j d\bar z_j$ be a $(0,1)$-form on $\Bn$. 
Then from (\ref{inverse}) the norm $|\alpha|_{i\partial\bar\partial\psi}$ of $\alpha$ with respect to $i\partial\bar\partial\psi$
can be written by
\begin{align*}
|\alpha|_{i\partial\bar\partial\psi}^2
&=\sum_{j,k=1}^n h_\psi^{jk}\alpha_j\bar\alpha_k\\
&=(1-|z|^2)^2\left(|\alpha|^2-\frac{2|\langle\alpha,  z\rangle|^2}{1+|z|^2}\right)\\
&=\frac{(1-|z|^2)^2}{1+|z|^2}\sum_{j,k=1}^n|\alpha_j\bar z_k-\alpha_k\bar z_j|^2
+\frac{(1-|z|^2)^3}{1+|z|^2}|\alpha|^2.
\end{align*}
Let $\beta=\sum_{j=1}^n\beta_j d z_j$ be a $(1,0)$-form on $\Bn$. 
Then the norm $|\beta|_{i\partial\bar\partial\psi}$ of $\beta$ with respect to $i\partial\bar\partial\psi$
can be written by
\begin{align*}
|\beta|_{i\partial\bar\partial\psi}^2
&=\sum_{j,k=1}^n\overline{h_\psi^{jk}}\beta_j\bar\beta_k\\
&=(1-|z|^2)^2\left(|\beta|^2-\frac{2|\langle\beta,  \bar z\rangle|^2}{1+|z|^2}\right)\\
&=\frac{(1-|z|^2)^2}{1+|z|^2}\sum_{j,k=1}^n|\beta_j z_k-\beta_k z_j|^2
+\frac{(1-|z|^2)^3}{1+|z|^2}|\beta|^2.
\end{align*}

We also obtain the following representation for the square root of $H_\psi(z)$,
\begin{align*}
H_\psi(z)^{1/2}=\frac{(1+|z|^2)^{1/2}}{(1-|z|^2)^{3/2}}\bar P_z + \frac{1}{1-|z|^2}\bar Q_z.
\end{align*}

\begin{lem}\label{Schwartz}
Let $\alpha$ be a $(0,1)$-form on $\Bn$.
Then
\begin{align}\label{Schwartz}
|\alpha|_{i\partial\bar\partial\psi}=\sup_{\xi\in\Cn\setminus\{0\}}\frac{|\langle\alpha, \xi\rangle|}{|\xi|_{h_\psi}}.
\end{align}
\end{lem}

\begin{proof}
All vectors in $\Cn$ are considered column vectors. Then
\begin{align*}
|\xi|_{h_\psi}^2=\xi^T H_\psi\bar\xi,\quad\xi\in\Cn.
\end{align*}
If we identify a $(0,1)$-form $\alpha=\sum_{j=1}^n\alpha_j d\bar z_j$ with the vector
$\alpha=(\alpha_1,\cdots, \alpha_n)^T$, then 
\begin{align*}
|\alpha|_{i\partial\bar\partial\psi}^2=\alpha^T H_\psi^{-1}\bar\alpha.
\end{align*}

We replace $\xi$ by $\overline{H_\psi(z)}^{-1/2}\,\xi$ to obtain 
\begin{align*}
\sup_{\xi\in\Cn\setminus\{0\}}\frac{|\langle\alpha, \xi\rangle|^2}{\xi^T H_\psi\bar\xi}
&=\sup_{\xi\in\Cn\setminus\{0\}}\frac{\big|\langle\alpha, \overline{H_\psi(z)}^{-1/2}\,\xi\rangle\big|^2}
{|\xi|^2}\\
&=\sup_{\xi\in\Cn\setminus\{0\}}\frac{\big|\langle\overline{H_\psi(z)}^{-1/2}\alpha, \xi\rangle\big|^2}
{|\xi|^2}\\
&=\big|\overline{H_\psi(z)}^{-1/2}\alpha\big|^2\\
&=\langle \overline{H_\psi(z)}^{-1}\alpha, \alpha\rangle\\
&=\alpha^T H_\psi^{-1}\bar\alpha.
\end{align*}
\end{proof}

The Hermitian metric $h_\psi$ defines a Riemannian metric $g_\psi$.
The metric $g_\psi$ is defined to be the real part of $h_\psi$:
\begin{align*}
g_\psi=\frac 12 (h_\psi +\bar h_\psi).
\end{align*}
Since $h_\psi(\xi,\xi) = g_\psi(\xi, \xi)$ for all $\xi\in\Cn$, the length
information of $h_\psi$ is entirely contained in its real part.

For a piecewise $C^1$ curve $\gamma : [0,1] \rightarrow \Bn$, the length induced by the Hermitian metric $h_\psi$ is defined by
\begin{align*}
\ell_{\psi}(\gamma):
= \int_{0}^{1} |\gamma'(t)|_{h_\psi}\, \,dt.
\end{align*} 
For $z, w \in \Bn$, the Riemannian distance induced by the Hermitian metric $h_\psi$ is
\begin{align*}
d_\psi(z,w):= \inf_\gamma\ell_{\psi}(\gamma) 
\end{align*} 
where $\gamma$ is a parametrized curve from $z$ to $w$ in $\Bn$.

\begin{pro}
Let $z\in\Bn$. Then
\begin{align*}
d_\psi(0,z)=\int_0^{|z|}\frac{(1+s^2)^{1/2}}{(1-s^2)^{3/2}}\,ds.
\end{align*}
\end{pro}

\begin{proof} 
Fix a point $z\in\Bn$.
Let $\gamma :[0,1]\rightarrow\Bn$ be a smooth and regular curve from $0$ to $z$.
According to part (d) of Lemma \ref{lem for H}, we have
\begin{align*}
\ell_{\psi}(\gamma) 
&= \int_{0}^{1} |\gamma'(t)|_{h_{\psi(\gamma(t))}}\,dt\\
&\ge\int_0^1\frac{(1+|\gamma(t)|^2)^{1/2}}{(1-|\gamma(t)|^2)^{3/2}}|P_{\gamma(t)}\gamma'(t)|\,dt.
\end{align*} 

Let $\alpha(t)=|\gamma(t)|$. Then $\alpha$ is smooth on $[0,1]$.
Since $\alpha^2(t)=\langle\gamma(t), \gamma(t)\rangle$, differentiation gives
\begin{align*}
2\alpha(t)\alpha'(t)=2{\rm Re}\langle\gamma'(t), \gamma(t)\rangle
=2{\rm Re}\langle P_{\gamma(t)}\gamma'(t), \gamma(t)\rangle.
\end{align*}
It follows that
\begin{align*}
|\alpha'(t)|\le|P_{\gamma(t)}\gamma'(t)|,\quad t\in[0,1].
\end{align*}
Hence
\begin{align*}
\ell_\psi(\gamma)
&\ge\int_0^1\frac{(1+|\gamma(t)|^2)^{1/2}}{(1-|\gamma(t)|^2)^{3/2}}|P_{\gamma(t)}\gamma'(t)|\,dt\\
&\ge\left|\int_0^1\frac{(1+\alpha^2(t))^{1/2}}{(1-\alpha^2(t))^{3/2}}\alpha'(t)\,dt\right|\\
&=\int_0^{|z|}\frac{(1+s^2)^{1/2}}{(1-s^2)^{3/2}}\,ds
\end{align*}
and so 
\begin{align*}
d_\psi(0,z)\ge\int_0^{|z|}\frac{(1+s^2)^{1/2}}{(1-s^2)^{3/2}}\,ds.
\end{align*}

We define $\gamma_0(t)=tz$ for $0\le t\le 1$. Then the length of $\gamma_0$ is
\begin{align*}
\ell_\psi(\gamma_0)
&=\int_0^1\frac{|z|(1+t^2|z|^2)^{1/2}}{(1-t^2|z|^2)^{3/2}}\,dt\\
&=\int_0^{|z|}\frac{(1+s^2)^{1/2}}{(1-s^2)^{3/2}}\,ds.
\end{align*}
This means that $\gamma_0$ is the geodesic curve connecting $0$ and $z$. 
So
\begin{align*}
d_\psi(0,z)=\int_0^{|z|}\frac{(1+s^2)^{1/2}}{(1-s^2)^{3/2}}\,ds.
\end{align*}
\end{proof}

\begin{cor}
Let $z\in\Bn$. Then
\begin{align*}
\frac{|z|}{\sqrt{1-|z|^2}}\le d_\psi(0,z)\le\frac{|z|(1+|z|^2)^{1/2}}{\sqrt{1-|z|^2}}.
\end{align*}
\end{cor}

\begin{pro}
The Hermitian metric space $(\Bn, h_\psi)$ is complete.
\end{pro}

\begin{proof}
Let 
\begin{align*}
f(z)=\log\left(\frac{1}{1-|z|^2}\right).
\end{align*}
Then $f$ is a smooth exhaustion function on $\Bn$.
Moreover, we have
\begin{align*}
\left|\bar\partial\log\left(\frac{1}{1-|z|^2}\right)\right|_{i\partial\bar\partial\psi}^2
=\frac{|z|^2(1-|z|^2)}{1+|z|^2}\le 1.
\end{align*}
That is, $f$ is a smooth exhaustion function with $|\bar\partial f|_{i\partial\bar\partial\psi}\le1$ on $\Bn$.
Let $z, w\in\Bn$ and $\gamma$ be a smooth curve in $\Bn$ connectiong $z$ and $w$.
By (\ref{Schwartz}), we have
\begin{align*}
|f(z)-f(w)|
&\le 2\int_0^1 |\langle\bar\partial f(\gamma(t)), \gamma'(t)\rangle|\, dt\\
&\le 2\int_0^1 |\bar\partial f|_{i\partial\bar\partial\psi}|\gamma'|_{h_{\psi(\gamma)}}\,dt\\
&\le 2\int_0^1 |\gamma'|_{h_{\psi(\gamma)}}\,dt.
\end{align*}
Hence
\begin{align*}
|f(z)-f(w)|\le 2 d_\psi(z,w).
\end{align*}
Since $f$ is an exhaustion function, all $d_\psi$-balls must be relatively compact in $\Bn$. 
Hence $(\Bn, h_\psi)$ is complete (see p. 366, (2.4) Lemma in \cite{Dem}).
\end{proof}


\subsection{The $d_\psi$- ball and polycylinder}
This section gives a geometric description of the $d_\psi$-ball $B_{\psi}(z, r)$ with polycylinder $D_{\psi}(z, r)$.

The $d_\psi$-ball centered at $z$ with radius $r>0$ is defined by the associated ball with $d_\psi$, namely, 
\begin{align*}
B_\psi(z,r) := \left\{w\in\Bn :  d_\psi(z,w) <r\right\}.
\end{align*}
We should note that $d_\psi$-ball is suitable tool for investigating exponential type weighted Bergman spaces on the unit ball rather than using the ball with radius function $\left(\Delta \psi\right)^{-\frac{1}{2}}$. Actually, the ball with $\left(\Delta \psi\right)^{-\frac{1}{2}}$ is helpful tool for studying function spaces with general exponential weights on ${\mathbf D}$ (\cite{Chr}, \cite{LR}, \cite{MO},
\cite{PP}) and ${\mathbf C^n}$ \cite{Dal}.
But it is not proper in the case of the unit ball.
For example, 
Lemma \ref{SMVP} with the reproducing property and comparable property implies the following estimate for the Bergman kernel on diagonal:
\begin{align}\label{Hessian}
K_{\psi}(z,z) \lesssim\frac{e^{\psi(z)}}{(1-\lvert z \rvert^2)^{2n+1}}.
\end{align}
The estimate is same as the result which can be obtained from Theorem 3.3 in \cite{CP} using series expansion. 
However, one could get only 
\begin{align}\label{Laplacian}
K_{\psi}(z,z) \lesssim\frac{e^{\psi(z)}}{(1-\lvert z \rvert^2)^{3n}}
\end{align} 
if one use the ball induced by the radius function with $\left(\Delta \psi\right)^{-\frac{1}{2}}$ instead of the $d_\psi$- ball 
$B_{\psi}(z,r)$. The estimate (\ref{Hessian}) is sharper than (\ref{Laplacian}) when $n>1$.  

We define the polycylinder $D_{\psi}(z,r)$ by
\begin{align*}
D_{\psi}(z,r)= \left\{w\in\Bn:\lvert z-P_z w\rvert < r\left(1-\lvert z \rvert^2\right)^{\frac32}, \, \lvert Q_z w\rvert < r\left(1-\lvert z \rvert^2\right) \right\}.
\end{align*}
When $n=1$, it follows that
$H_\psi (z) = \Delta \psi (z) \simeq (1-\lvert z \rvert^2)^{-3}$.
Moreover, note that $P_z = I$ for $z\in{\D}\setminus\{0\}$, 
and $Q_z = 0$.
Hence the polycylinder $D_{\psi}(z,r)$, the $d_\psi$-ball $B_\psi(z,r)$, and
the ball with radius function $\left(\Delta \psi\right)^{-\frac{1}{2}}$ are all same when $n=1$.

Let $B(z,r)=\{ w\in\Bn : |w-z|<r\}$ be the Euclidean ball centered at $z$ of radius $r$.

\begin{lem}\label{CarlesonBox}
Let $z\in\Bn$ and $0<r<\frac 12$. Then
\begin{align*}
(1-2r)(1-|z|^2)\le 1-|w|^2\le(1+2r)(1-|z|^2)
\end{align*}
for every $w\in B(z, r(1-|z|^2))$.
\end{lem}

\begin{cor}\label{comp in D}
Let $0< r <\frac{1}{4}$.
Then
\begin{align*}
(1-4r)(1-|z|^2)\le 1-|w|^2\le(1+4r)(1-|z|^2)
\end{align*} 
whenever $w \in D_{\psi}(z,r)$ for $z, w \in \Bn$.
\end{cor}

\begin{proof}
Let $w \in D_{\psi}(z,r)$. Then
$\lvert z- w\rvert \le \lvert z-P_z w\rvert+\lvert Q_z w\rvert < 2r(1-\lvert z \rvert^2)$
which means $w$ belongs $B\left(z, 2r(1-\lvert z \rvert^2)\right)$. 
It implies 
\begin{align*}
(1-4r)(1-|z|^2)\le 1-|w|^2\le(1+4r)(1-|z|^2)
\end{align*} 
provided $0<2r < \frac 12$ by Lemma \ref{CarlesonBox}.
\end{proof}

\begin{pro}\label{comparable1}
Let $0< r < \frac 54$. Then
\begin{align*}
D_{\psi}\left(z, \frac{r}{10}\right) \subset B_\psi(z,r).
\end{align*}
\end{pro}

\begin{proof}
Recall that the distance induced by the metric $h_{\psi}$ is   
\begin{align*}
d_\psi(z,w) = \inf_\gamma \ell_{\psi}(\gamma)
\end{align*} 
where the infimum is taken over all piecewise smooth curve ${\gamma}:[0,1] \rightarrow \Bn$ with $\gamma(0)=z$ and $\gamma(1)=w$ and the length induced by the Hessian metric $h_\psi$ is 
\begin{align*}
\ell_{\psi}(\gamma) = \int_{0}^{1} \left\{ \frac{2\left|\langle \gamma(t),\gamma'(t)\rangle\right|^2}
{\left(1-|\gamma(t)|^2\right)^{3}}+\frac{\left|\gamma'(t)\right|^2}{\left(1-|\gamma(t)|^2\right)^{2}} \right\}^{\frac{1}{2}} dt
\end{align*} 
for each curve $\gamma$.

Suppose that $w$ belongs to $D_{\psi}(z,m)$ with $0<m<1$. We assume $w \neq z$ without loss of generality.
We have 
\begin{align*}
(1-4m)(1-|z|^2)\le 1-|w|^2\le(1+4m)(1-|z|^2)
\end{align*} 
when $m<\frac 14$ by Lemma \ref{comp in D}. 
Let $\widehat{\gamma}_1$ be a line segment from $z$ to $P_z w$ and $\widehat{\gamma}_2$ be a line segment from $P_z w$ to $w$, 
precisely, 
\begin{align*}
\widehat{\gamma}_1 (t)=tP_z w+(1-t)z
\end{align*} 
and 
\begin{align*}
\widehat{\gamma}_2 (t)=tw+(1-t)P_z w.
\end{align*} 
Let $\widehat{\gamma}$ be a parametrized curve for $\widehat{\gamma}_1 + \widehat{\gamma}_2$. Then
$\widehat{\gamma}$ is a curve connecting $z$ and $w$ in $D_\psi(z,r)$ and 
\begin{align*}
\ell_{\psi}(\widehat{\gamma})=\ell_{\psi}(\widehat{\gamma}_1)+\ell_{\psi}(\widehat{\gamma}_2).
\end{align*} 

We have
\begin{align*}
\ell_{\psi}(\widehat{\gamma}_1) 
&\le  \int_{0}^{1} \left\{\frac{\sqrt{2}}{\left(1-|\widehat{\gamma}_1(t)|^2\right)^{\frac{3}{2}}}\left|\langle 
\widehat{\gamma}_1(t),\widehat{\gamma}_1'(t)\rangle\right| 
+ \frac{1}{1-|\widehat{\gamma}_1(t)|^2} \left|\widehat{\gamma}_1'(t)\right|\right\} dt \\ 
&\le \frac{4}{\left(1-|z|^2\right)^{\frac{3}{2}}} \int_{0}^{1} \left|\langle \widehat{\gamma}_1(t),\widehat{\gamma}_1'(t)\rangle\right| dt + \frac{2}{1-|z|^2} \int_{0}^{1} \left|\widehat{\gamma}_1'(t)\right| dt 
\end{align*}
for $0<m<\frac 18$ by Lemma \ref{comp in D}. 
Note that $\widehat{\gamma}_1'(t) = P_z w - z $. The Cauchy-Schwartz inequality yields that 
\begin{align*}
\ell_{\psi}(\widehat{\gamma}_1) 
\le& \frac{4}{\left(1-|z|^2\right)^{\frac{3}{2}}} \int_{0}^{1} \left|\widehat{\gamma}_1(t)\right| \left|\widehat{\gamma}_1'(t)\right| dt + 2\frac{|z-P_z w|}{1-|z|^2} \\ 
\le& \frac{4}{\left(1-|z|^2\right)^{\frac{3}{2}}} \sup_{t} \left|\widehat{\gamma}_1(t)\right| \int_{0}^{1} \left|\widehat{\gamma}'_1(t)\right| dt + 2\frac{|z-P_z w|}{1-|z|^2}. 
\end{align*}
The fact $\sup_{t} \left|\widehat{\gamma}_1(t)\right| \le 1$ gives
\begin{align*}
\ell_{\psi}(\widehat{\gamma}_1) \le \frac{4}{\left(1-|z|^2\right)^{\frac{3}{2}}} |z-P_z w|+ 2\frac{|z-P_z w|}{1-|z|^2}
< 6m
\end{align*}
when $w \in D_{\psi}(z,m).$

We also have 
\begin{align*}
\ell_{\psi}(\widehat{\gamma}_2)
&\le \frac{4}{\left(1-|z|^2\right)^{\frac{3}{2}}} \int_{0}^{1} \left|\langle \widehat{\gamma}_2(t),\widehat{\gamma}_2'(t)\rangle\right| dt + \frac{2}{1-|z|^2} \int_{0}^{1} \left|\widehat{\gamma}_2'(t)\right| dt
\end{align*}
by Lemma \ref{comp in D}. Note that $\widehat{\gamma}_2'(t) = Q_z w$.
The fact that $P_z w$ and $Q_z w$ are perpendicular asserts
\begin{align*}
\langle \widehat{\gamma}_2(t),\widehat{\gamma}_2'(t)\rangle = \left\langle (Q_z w)t+P_z w, Q_z w\right\rangle = t \left|Q_z w\right|^2.
\end{align*}
Then, 
\begin{align*}
\ell_{\psi}(\widehat{\gamma}_2) 
\le& \frac{4}{\left(1-|z|^2\right)^{\frac{3}{2}}} |Q_z w|^2 \int_{0}^{1} t \,dt + 2\frac{|Q_z w|}{1-|z|^2}\\
\le& \frac{2}{\left(1-|z|^2\right)^{\frac{3}{2}}} |Q_z w|^2 + 2\frac{|Q_z w|}{1-|z|^2} \\
<& 4m
\end{align*} 
when $w \in D_{\psi}(z,m).$

As a result, we get 
\begin{align*}
d_\psi(z,w) \le \ell_{\psi}(\widehat{\gamma}) < 10m
\end{align*} 
which implies $w\in  B_{\psi}(z, 10m)$. By putting $m=\frac{r}{10}$, it is obtained   
\begin{align*}
D_{\psi}\left(z, \frac{r}{10}\right) \subset B_\psi(z,r).
\end{align*}
\end{proof}

\begin{pro}\label{ball_size_comp2}
Let $0< r<\frac{1}{12}$. Then
\begin{align*}
B_\psi(z,r) \subset D_{\psi}(z, 2r).
\end{align*}
\end{pro}

\begin{proof}
Suppose that $w$ belongs to $B_\psi(z,r)$. We assume $w \neq z$ without loss of generality. The proof is divided into three steps.

{\bf Step 1.} We will show that $d_\psi(z, w) <r$ implies $|z-w| < \frac 32 r(1-|z|^2)$ and $|Q_z(w)| \le \frac 32 r(1-|z|^2)$. 

Suppose $d_\psi(z, w) <r$. 
For any piecewise $C^1$ curve $\gamma$ from $z$ to $w$, let $T_0 \in (0,1]$ be the minimum of $t$ satisfying 
\begin{align*}
|z- \gamma(t)| = \min\left\{|z-w|, s(1-|z|^2)\right\},
\end{align*}
where $s$ is an undetermined positive real number.
Then 
\begin{align*}
|z- \gamma(t)|\le |z- \gamma(T_0)| \le s(1-|z|^2)\quad \text{for}\quad t \in [0,T_0].
\end{align*}
By Lemma \ref{CarlesonBox}, if $0<s<\frac 12$, then it gives
\begin{align*}
(1-2s)(1-|z|^2) \le 1-|\gamma(t)|^2 \le (1+2s)(1-|z|^2) \quad \text{for}\quad t \in [0,T_0].
\end{align*} 
We have
\begin{align}\label{Euclid}
\begin{split}
\ell_{\psi}(\gamma) 
&\ge \int_{0}^{1} \frac{1}{1-|\gamma(t)|^2} \left|\gamma'(t)\right| dt\\
&\ge \int_{0}^{T_0} \frac{1}{1-|\gamma(t)|^2} \left|\gamma'(t)\right| dt\\
&\ge \frac{1}{(1+2s)(1-|z|^2)}|z-\gamma(T_0)|.
\end{split}
\end{align}
Since $\gamma$ is an arbitrary curve connecting $z$ and $w$, 
\begin{align*}
d_\psi(z,w) \ge \frac{1}{1+2s} \min\left\{\frac{|z-w|}{1-|z|^2}, s\right\}.
\end{align*}
We take $s=\frac 14$. Then we have
\begin{align*}
r > d_\psi(z,w) \ge \frac{1}{1+2s} \min\left\{\frac 23\frac{|z-w|}{1-|z|^2}, \frac 16\right\},\quad w \in B_\psi(z,r).
\end{align*}
Since $0<r<\frac 16$, it is contradiction whenever 
\begin{align*}
\frac 23\frac{|z-w|}{1-|z|^2}\ge\frac 16.
\end{align*}
Hence 
\begin{align*}
\frac 23\frac{|z-w|}{1-|z|^2}<\frac 16
\end{align*}
and we have 
\begin{align*}
r > d_\psi(z,w) \ge \frac 23\frac{|z-w|}{1-|z|^2}.
\end{align*}
It gives 
\begin{align}\label{B_H radius}
|z-w| < \frac 32 r(1-|z|^2) \quad \text{for} \quad w \in B_\psi(z,r).
\end{align}
Moreover, 
\begin{align}\label{concl_step1}
|Q_z w| \le |z-w|<\frac 32 r(1-|z|^2) \quad \text{for} \quad w \in B_\psi(z,r)
\end{align}
since $z-P_z w$, $z-w$, and $Q_z w$ construct a right triangle with hypotenuse $z-w$.

Now, we will show that $d_\psi(z, w) <r$ implies $|z-P_z w| \lesssim r(1-|z|^2)^{\frac{3}{2}}$ when $|z| \le \frac{1}{2}$ and $|z| > \frac{1}{2}$ in Step 2 and Step 3, respectively. 

{\bf Step 2.} 
We assume that $|z| \le \frac{1}{2}$. Then $\frac 34\le 1-|z|^2\le 1$ and
\begin{align*}
|z-w|<\frac 32r(1-|z|^2)\le\sqrt{3}r(1-|z|^2)^{3/2}
\end{align*}
by \eqref{B_H radius} in Step 1. 
Since $z-P_z w$, $z-w$, and $Q_z w$ construct a right triangle with hypotenuse $z-w$, we have 
\begin{align}\label{concl_step2}
|z-P_z w| <\sqrt{3}r(1-|z|^2)^{3/2}  \quad \text{for} \quad w \in B_\psi(z,r)
\end{align}
when $|z| \le \frac{1}{2}$.
 
{\bf Step 3.} 
We assume that $|z| > \frac{1}{2}$.
Now suppose that $0<r<\frac {1}{12}$ and $d_\psi(z, w)<r$. 
We also assume $z \neq P_z w$ without loss of generality. 
Hereinafter, we consider only the curves $\gamma$ connecting $z$ and $w$ satisfying 
$$
\ell_{\psi}(\gamma) \le 2 d_\psi(z, w).
$$ 
This means that
\begin{align*}
\gamma(t)\in B_\psi(z, 2r)\quad\text{for every}\quad t.
\end{align*}

For each curve $\gamma$, define $\gamma_1(t)=P_z(\gamma(t))$ and $\gamma_2(t)=\gamma(t)-\gamma_1(t)= Q_z(\gamma(t))$.
Let $t_0 \in [0,1]$ be the minimum of $t$ such that 
\begin{align*}
|z- \gamma_1(t)| = |z-P_z w|. 
\end{align*}
Since $0<2r<\frac 16$, by \eqref{B_H radius} and Lemma \ref{CarlesonBox}, we have
\begin{align}\label{comp in B_H}
\frac 12(1-|z|^2)\le 1-|w|^2\le\frac 32(1-|z|^2)
\end{align} 
for $w \in B_\psi(z, 2r)$. 
It gives
\begin{align*}
\ell_{\psi}(\gamma) 
&\ge \int_{0}^{1} \frac{\sqrt{2}}{\left(1-|\gamma(t)|^2\right)^{\frac{3}{2}}}\left|
\langle \gamma(t),\gamma'(t)\rangle\right| \,dt\\
&\ge\frac{2\sqrt{2}}{3}\frac{1}{\left(1-|z|^2\right)^{\frac{3}{2}}} 
\int_{0}^{t_0} \left|\langle \gamma(t),\gamma'(t)\rangle\right|\, dt\\
&\ge\frac{2\sqrt{2}}{3} \frac{1}{\left(1-|z|^2\right)^{\frac{3}{2}}} 
\left|\int_{0}^{t_0}|\langle \gamma_1(t),\gamma_1'(t)\rangle|-|\langle \gamma_2(t),\gamma_2'(t)\rangle|\,dt\right|.
\end{align*}

Since $\gamma_1$ and $\gamma_1'$ are parallel, 
\begin{align*}
\int_{0}^{t_0} \left|\langle \gamma_1(t),\gamma_1'(t)\rangle\right| dt =\int_{0}^{t_0} \left|\gamma_1(t)\right|\left|\gamma_1'(t)\right| dt \ge \inf_{0\le t \le t_0} \left|\gamma_1(t)\right| \int_{0}^{t_0} \left|\gamma_1'(t)\right| dt. 
\end{align*}
The hypothesis $\ell_{\psi}(\gamma) \le 2 d_\psi(z, w)$ gives 
\begin{align*}
d_\psi(z, \gamma_1 (t)) < d_\psi(z, \gamma (t)) < 2r.
\end{align*}
Inequality \eqref{B_H radius} and $0 < 1-|z|^2 < \frac{3}{4}$ yield
\begin{align*}
\gamma_1 (t) \in B_\psi(z, 2r) \subset B (z, 3r(1-|z|^2)) \subset B \left(z, \frac 94 r\right).
\end{align*}
Hence, we have 
\begin{align*}
\inf_{0\le t \le t_0} \left|\gamma_1(t)\right| \ge |z| - \frac 94 r > \frac {5}{16}
\end{align*}
with $|z|>\frac 12$ and $0<r<\frac {1}{12}$.
It is obvious that $\int_{0}^{t_0} \left|\gamma_1'(t)\right| dt \ge |z-P_z w|$.
We have 
\begin{align*}
\int_{0}^{t_0} \left|\langle \gamma_1(t),\gamma_1'(t)\rangle\right| dt 
\ge \frac{5}{16} |z-P_z w|. 
\end{align*}

By the Cauchy-Schwartz inequality,
\begin{align*}
\int_{0}^{t_0} \left|\langle \gamma_2(t),\gamma_2'(t)\rangle\right| dt
\le \int_{0}^{t_0} \left| \gamma_2(t) \right| \left| \gamma_2'(t)\right| dt
\le  \sup_{0\le t \le t_0} \left| \gamma_2(t) \right| \int_{0}^{t_0} \left| \gamma_2'(t)\right| dt.
\end{align*}
Since $z-\gamma_1(t)$, $z-\gamma(t)$ and $\gamma_2(t)$ 
construct a right triangle with hypotenuse $z-\gamma(t)$, we have
\begin{align*}
\sup_{0\le t \le t_0} \left| \gamma_2(t) \right|\le \sup_{0\le t \le t_0} |z-\gamma(t)| \le |z-\gamma(t_0)| 
\le \int_{0}^{t_0} \left| \gamma'(t)\right| dt. 
\end{align*}
Thus we get
\begin{align*}
\int_{0}^{t_0} \left|\langle \gamma_2(t),\gamma_2'(t)\rangle\right| dt
\le \left\{\int_{0}^{t_0} \left| \gamma'(t)\right| dt\right\}^2
\le \frac 94(1-|z|^2)^2 \ell_{\psi}(\gamma)^2
\end{align*}
since 
\begin{align*}
\int_{0}^{t_0} \left|\gamma'(t)\right| dt \le \frac 32(1-|z|^2) \ell_{\psi}(\gamma) 
\end{align*}
as in \eqref{Euclid} with (\ref{comp in B_H}).

Therefore, we obtain  
\begin{align*}
\ell_{\psi}(\gamma)
&\ge \frac{5\sqrt{2}}{24} \frac{1}{\left(1-|z|^2\right)^{\frac{3}{2}}} |z-P_z w| - \frac{3\sqrt{2}}{2}\left(1-|z|^2\right)^{\frac{1}{2}} \ell_{\psi}(\gamma)^2
\end{align*}
which allows
\begin{align*}
\frac{5\sqrt{2}}{24} \frac{1}{\left(1-|z|^2\right)^{\frac{3}{2}}} |z-P_z w|
&\le \ell_{\psi}(\gamma)+\frac{3\sqrt{2}}{2}(1-|z|^2)^{1/2}\ell_\psi(\gamma)^2\\
&\le \frac 12 d_\psi(z,w)+\frac{3\sqrt{2}}{2}(1-|z|^2)^{1/2}\frac 14 d_\psi(z,w)^2\\
&\le \frac 12 r+\frac{3\sqrt{2}}{8}(1-|z|^2)^{1/2} r^2.
\end{align*}
It shows
\begin{align}
\begin{split}\label{concl_step3}
|z-P_z w| &< \frac{12}{5\sqrt{2}}r \left(1-|z|^2\right)^{\frac{3}{2}} 
+\frac 95 r^2 (1-|z|^2)^2\\
&\le\left(\frac{12}{5\sqrt{2}}+\frac{9}{5}r\right)r(1-|z|^2)^{\frac 32}\\
&\le2r(1-|z|^2)^{\frac 32}\quad \text{for} \quad w \in B_\psi(z,r)
\end{split}
\end{align}
when $0<r<\frac{1}{12}$.

Finally, we get the desired result 
$$
B_\psi(z,r) \subset D_{\psi}(z, 2r) 
$$ 
by gathering with \eqref{concl_step1}, \eqref{concl_step2}, and \eqref{concl_step3}. 
\end{proof}

From Propositions \ref{comparable1} and \ref{ball_size_comp2}, we have the following theorem.

\begin{thm}\label{ball_size_comp}
Let $0< r <\frac{1}{12}$. Then 
\begin{align*}
D_{\psi}\left(z, \frac{r}{10}\right) \subset B_\psi(z,r) \subset D_{\psi}(z, 2r).
\end{align*}
\end{thm}

For non zero $z\in\C^n$, let 
\begin{align*}
U_1=\frac{\bar z}{|z|}.
\end{align*}
Let $U_j\in\C^n$ be such that
\begin{align*}
\langle z, \bar U_j\rangle=0,\quad j\ge 2
\end{align*}
and
\begin{align*}
\langle U_j, U_k\rangle=\delta_{jk}
\quad 2\le j, k \le n.
\end{align*}
Then 
\begin{align*}
U_z=(U_1, \cdots, U_n)^T
\end{align*}
is a unitary matrix such that
\begin{align*}
U_zz=(|z|, 0, \cdots, 0).
\end{align*}
We note that
\begin{align*}
U_z(P_zw)=\left(\frac{\langle w, z\rangle}{|z|}, 0, \cdots, 0\right)
\end{align*}
and 
\begin{align*}
U_z(Q_zw)=U_zw-\left(\frac{\langle w, z\rangle}{|z|}, 0, \cdots, 0\right).
\end{align*}
Let $\zeta=U_zw$. Then
\begin{align*}
\zeta_1=\frac{\langle w, z\rangle}{|z|}
\end{align*}
and
\begin{align*}
U_z(D_{\psi}(z,r))
=\left\{\zeta :\left||z|-\zeta_1\right| < r\left(1-\lvert z \rvert^2\right)^{\frac32}, \, \lvert\zeta'\rvert < r\left(1-\lvert z \rvert^2\right) \right\},
\end{align*}
where $\zeta'=(\zeta_2,\cdots,\zeta_n)$.
We note that $U_z(D_{\psi}(z,r))$ is a polycylinder 
which is a topological product of a disk in $\mathbf C$ and
a ball in $\mathbf C^{n-1}$.

\begin{cor}\label{cor_vol}
For $0< r_1, r_2 <\frac {1}{12}$,
\begin{align}\label{comp volume}
{\rm Vol}(B_{\psi}(z, r_1)) \simeq {\rm Vol}(D_{\psi}(z,r_2)) \simeq \left(1-|z|^2\right)^{2n+1}.
\end{align}
\end{cor}

\subsection{Test functions}

For $z \in \Bn$, the involutive automorphisms  on $\Bn$ are defined 
\begin{align*}
 \varphi_{z}(w):=\frac{z-P_z w -\sqrt{1-\lvert z \rvert^2}Q_z w}
{1-\langle w, z \rangle}.   
\end{align*} 
It has the following property:
\begin{align}\label{prop4_phi}
  1-\lvert \varphi_{z} (w)\rvert^2 = \frac{(1-\lvert z \rvert^2)(1-\lvert w \rvert^2)}{\lvert 1-\langle w, z\rangle \rvert^2}.
\end{align}
For more details of the automorphisms of $\Bn$, see page 23 of \cite{Zhu}. 

Due to the definition of $D_{\psi}(z,r)$, we can get the following inequality which is essential for proving the estimate of the test function in Lemma \ref{test fn lem}. 

\begin{lem}\label{imp_ineq}
For $z \in \Bn$ and small $r > 0$, there is a constant $C$ depending only on the radius $r$ satisfying
\begin{align*}
\left| 2 \text{\rm Re} \left(\frac{1}{1- \langle w, z \rangle }\right) -\frac{1}{1-\lvert z \rvert^2} -\frac{1}{1-\lvert w \rvert^2}\right|\le C
\end{align*}
for $w \in D_{\psi}(z, r)$. 
\end{lem}

\begin{proof}
Using \eqref{prop4_phi}, we get the reformulation:
\begin{align}\label{eq_difference}
2\re &\left(\frac{1}{1- \langle w, z \rangle }\right) -\frac{1}{1-\lvert z \rvert^2} -\frac{1}{1-\lvert w \rvert^2} \nonumber\\
& = \frac{\lvert z-w\rvert^2}{\lvert 1- \langle w ,z\rangle \rvert^2} - \lvert \varphi_z (w)\rvert^2 \left(\frac{1}{1-\lvert z \rvert^2} + \frac{1}{1-\lvert w \rvert^2}\right),
\end{align}
which indicates 
\begin{align}\label{Eq}
- \lvert \varphi_z (w)\rvert^2 \left( \frac{1}{1-\lvert z \rvert^2} + \frac{1}{1-\lvert w \rvert^2}\right) \le \text{LHS of \eqref{eq_difference}} \le \frac{\lvert z-w\rvert^2}{\lvert 1- \langle w ,z\rangle \rvert^2}.
\end{align}

First, we show that $\frac{\lvert z-w\rvert^2}{\lvert 1- \langle w ,z\rangle \rvert^2}$ is dominated with some constant independent of $z$ and $w$, 
which means the LHS of \eqref{eq_difference} has an upper bound $C_r$. 
For $w\in D_{\psi}(z,r)$, we have $\lvert z-w\rvert^2 < 4r^2(1-\lvert z \rvert^2)^2$ since $\lvert z-w\rvert \le \lvert z-P_z (w)\rvert+\lvert Q_z w\rvert$.

By Lemma \ref{comp in D}, we have $1-\lvert z \rvert^2 \simeq 1-\lvert w \rvert^2$ for $w \in D_{\psi}(z,r)$ for small $r>0$.
Hence there exists $C_r>0$ such that
\begin{align}\label{upper constant 1}
\frac{\lvert z-w\rvert^2}{\lvert 1- \langle w ,z\rangle \rvert^2} < \frac{4 r^2(1-\lvert z \rvert^2)^2}{\lvert 1- \langle w ,z\rangle \rvert^2} \le C_r \frac{(1-\lvert z \rvert^2)(1-\lvert w \rvert^2)}{\lvert 1- \langle w ,z\rangle \rvert^2}.
\end{align}
The RHS of \eqref{upper constant 1} is equal to $$C_r \left(1-\lvert \varphi_{z} (w)\rvert^2\right)$$ by \eqref{prop4_phi}. It is dominated by $C_r$ since $\varphi_{z} (w)$ belongs to $\Bn$. 

Next, we show that the LHS of \eqref{Eq} has a lower bound $- C_r'$. Since $z-P_z (w)$ and $Q_z (w)$ are perpendicular, we have 
\begin{align*}
\lvert \varphi_z(w)\rvert^2 = \frac{\lvert z-P_z (w)\rvert^2 + (1-\lvert z \rvert^2) 
\lvert Q_z (w)\rvert^2}{\lvert 1- \langle w ,z\rangle \rvert^2}. 
\end{align*}
The definition of $D_{\psi}(z, r)$ yields
$$\lvert \varphi_z (w)\rvert^2 < \frac{2r^2(1-\lvert z \rvert^2)^3}{\lvert 1- \langle w ,z\rangle \rvert^2}$$ 
for $w \in D_{\psi}(z, r)$.
It implies 
\begin{align*}
\lvert \varphi_z (w)\rvert^2 \left( \frac{1}{1-\lvert z \rvert^2} + \frac{1}{1-\lvert w \rvert^2}\right)
&< \frac{2r^2(1-\lvert z \rvert^2)^3}{\lvert 1- \langle w ,z\rangle \rvert^2}\left( \frac{1}{1-\lvert z \rvert^2} + \frac{1}{1-\lvert w \rvert^2}\right)\\
&\le C_r' \frac{(1-\lvert z \rvert^2)(1-\lvert w \rvert^2)}{\lvert 1- \langle w ,z\rangle \rvert^2} \le C_r'
\end{align*}
with Lemma \ref{comp in D}.
The proof is done by getting $C=\max \left\{C_r , C_r'\right\}$.
\end{proof}

By the previous inequality, we can get the following lemma for test functions. It will be used for proving the weighted sub-mean-value property and the results in Section 3 and Section 4.  

\begin{lem}\label{test fn lem}
For $z \in \Bn$, let $F_{z} (w) := e^{\frac{1}{1-\langle w,z\rangle}-\frac{1}{2}\frac{1}{1-\lvert z \rvert^2}}$. The holomorphic function $F_{z}$ has following properties:
\begin{align}\label{test fn_1}
\lvert F_{z} (w)\rvert^2 e^{-\frac{1}{1-\lvert w \rvert^2}} \simeq 1 \quad \text{when} \quad w \in D_{\psi}(z, r). 
\end{align}
\end{lem}

\begin{proof}
By Lemma \ref{imp_ineq}, we can show that for $z \in \Bn$ and small $r > 0$, there is a constant $C$ depending only $r$ satisfying
  \begin{align*}
{C}^{-1} e^{- \frac{1}{1-\lvert z \rvert^2} - \frac{1}{1-\lvert w \rvert^2}} \le \lvert  e^{-\frac{1}{1- \langle w, z \rangle }}\rvert^{2} \le C e^{- \frac{1}{1-\lvert z \rvert^2} - \frac{1}{1-\lvert w \rvert^2}}
\end{align*}
for $w \in D_{\psi}(z, r)$. It gives \eqref{test fn_1}.
\end{proof}


\subsection{Sub-mean-value inequality}

\begin{lem}\label{SMVP}
Let $f$ be a holomorphic function on $\Bn$ 
For $z\in\Bn$ and a sufficiently small radius $r>0$, there is a constant $C$ depending on $r$ satisfying 
\begin{align*}
\lvert f(z)\rvert^2\le\frac{C}{(1-\lvert z \rvert^2)^{2n+1}} 
\int_{D_{\psi}(z,r)} \lvert f(\zeta)\rvert^2\,dV(\zeta). 
\end{align*}
\end{lem}

\begin{proof}
We note that
\begin{align*}
\lvert f(z)\rvert^2 
&=\lvert f\circ U_z^{-1}(|z|, 0, \cdots, 0)\rvert^2\\
&\leq\frac{C}{{\rm Vol}(D_\psi((|z|,\cdots,0),r))} \int_{D_\psi((|z|,\cdots,0),r)}\lvert f\circ U_z^{-1}(\zeta)\rvert^2\,dV(\zeta)\\
&=\frac{C}{{\rm Vol}(D_{\psi}(z,r))}\int_{D_{\psi}(z,r)} \lvert f(\zeta)\rvert^2\,dV(\zeta).
\end{align*}
\end{proof}

\begin{lem}\label{submean}
Let $f$ be a holomorphic function on $\Bn$ and $s \in {\R}$. 
For $z \in \Bn$ and a sufficiently small radius $r>0$, there is a constant $C$ depending on $s$ and $r$ satisfying 
\begin{align*}
\lvert f(z)\rvert^2 e^{-\frac{s}{1-\lvert z \rvert^2}} \le \frac{C}{(1-\lvert z \rvert^2)^{2n+1}} \int_{B_{\psi}(z,r)} \lvert f(\zeta)\rvert^2 e^{-\frac{s}{1-\lvert \zeta\rvert^2}} \,dV(\zeta).  
\end{align*}
\end{lem}

\begin{proof}
Since the function $F_{z}$ is non-vanishing, $f(\zeta)F_{z}(\zeta)^{-s}$ with a principle branch is holomorphic, and 
$\lvert f(\zeta)F_{z}(\zeta)^{-s} \rvert^2$ is plurisubharmonic. From Lemma \ref{SMVP} there is $\delta > 1$ satisfying
\begin{align*}
\lvert f(z) F_{z}(z)^{-s}\rvert^2 
&\leq \frac{C}{{\rm Vol}(D_{\psi}(z,\delta^{-1}r))} \int_{D_{\psi}(z,\delta^{-1}r)} \lvert f(\zeta)\rvert^2 \lvert F_{z}(\zeta)^{-2}\rvert^{s} \,dV(\zeta)\\
&\le \frac{C}{{\rm Vol}(D_{\psi}(z,\delta^{-1}r))} \int_{D_{\psi}(z,\delta^{-1}r)} \lvert f(\zeta)\rvert^2 e^{-\frac{s}{1-\lvert \zeta\rvert^2}} \,dV(\zeta) 
\end{align*}
with aid of \eqref{test fn_1}.
We note that  $\lvert f(z)\rvert^2 e^{-\frac{s}{1-\lvert z \rvert^2}} = \lvert f(z) F_{z}(z)^{-s}\rvert^2$.
Theorem \ref{ball_size_comp} and Corollary \ref{cor_vol} yield
\begin{align*}
\lvert f(z)\rvert^2 e^{-\frac{s}{1-\lvert z \rvert^2}} 
&\le \frac{C}{(1-\lvert z \rvert^2)^{2n+1}} \int_{B_{\psi}(z,r)} \lvert f(\zeta)\rvert^2 e^{-\frac{s}{1-\lvert \zeta\rvert^2}} \,dV(\zeta).
\end{align*}
\end{proof}

\section{Estimates for the Bergman kernel}

\subsection{Auxiliary lemmas}

Let ${\rm Lip}(\Bn, h_\psi)$ be the class of scalar functions
$u : \Bn\rightarrow\R$ that are Lipschitz continuous with respect to $d_\psi$.
Let
\begin{align*}
{\rm Lip}_\psi(u)=\sup\left\{\frac{|u(z)-u(w)|}{d_\psi(z,w)} : z, w\in\Bn, z\neq w\right\}.
\end{align*}
Then by Rademacher's theorem, $u$ is almost everywhere differentiable and, by our convention
$g_\psi={\rm Re}\, h_\psi$, 
\begin{align*}
\sup_{\Bn}|\partial u|_{i\partial\bar\partial\psi}=\sup_{\Bn}|\bar\partial u|_{i\partial\bar\partial\psi}\le\frac 12 {\rm Lip}_\psi(u).
\end{align*}

\begin{lem}\label{eta}
There exists $\eta\in C^\infty_c(\R, \R)$ such that
\begin{itemize}
\item[(1)]
$|\eta(x)|\le 1,\quad x\in\R$.
\item[(2)]
$\eta(x)=1,\quad |x|\le 1$.
\item[(3)]
${\supp}(\eta)\subset (-2, 2)$.
\item[(4)]
$(\eta'(x))^2\le C \eta(x),\quad x\in\R$.
\end{itemize}
\end{lem}

\begin{proof}
We define
\begin{align*}
f(x)=
\begin{cases}
e^{-\frac 1x},\quad x>0\\
0,\quad x\le 0
\end{cases}
\end{align*}
and
\begin{align*}
\eta(x)=\frac{f(2-|x|)}{f(|x|-1)+f(2-|x|)}.
\end{align*}
Then $\eta\in C^\infty_c(\R, \R)$ satisfies the required conditions.
\end{proof}

\begin{lem}\label{char_fn_2}
For a small $\delta$, there is a characteristic function $\chi_z : \C^n \rightarrow \R$ depending on $\varphi$ with the following properties:
\begin{itemize}
\item [(a)] $0\leq \chi_z \leq1$.
\item [(b)] $\chi_z |_{B_{\psi}(z, \frac{\delta}{4})} \equiv 1$.
\item [(c)] $\supp \chi_z \subset B_{\psi}(z, \delta)$.
\item [(d)] $|\bar\partial {\chi_z}|_{i\partial\bar\partial\psi}^2 \leq C \frac{1}{\delta^2} \chi_z$ for some $C>0$. 
\end{itemize}
\end{lem}

\begin{proof}
We note that $(\Bn, h_\psi)$ is a complete Hermitian manifold and $d_\psi(z, \cdot)\in {\rm Lip}(\Bn, h_\psi)$ with
\begin{align*}
{\rm Lip}_{\psi}(d_\psi(z, \cdot))=1
\end{align*}
by the triangular inequality.
Let $\kappa:(\Bn, h_\psi) \rightarrow (0, +\infty)$ be a continuous function.
By a reformulated Greene-Wu's theorem (see \cite{GW} and \cite{AFLR}), 
there is 
a smooth real-valued Lipschitz function $u_z$ 
satisfying  
	\begin{align*}
	\left|\partial u_z\right|_{i\partial\bar\partial\psi}=\left|\bar\partial u_z\right|_{i\partial\bar\partial\psi} \le 1
	\end{align*} 
	and 
	\begin{align}\label{kerBn_smooth_fn_1}
	\left|u_z - d_\psi(z, \cdot)\right| \le\frac{\kappa(\cdot)}{2}. 
	\end{align}
Let $\chi_z(\zeta) = \eta(\frac{4}{3\delta} u_z(\zeta))$, where $\eta$ is the function defined in Lemma \ref{eta}.
Then
\begin{align*}
\bar\partial\chi_z(\zeta)=\frac{4}{3\delta}\eta'\left(\frac{4}{3\delta} u_z(\zeta)\right)\bar\partial u_z(\zeta)
\end{align*}
and so
\begin{align*}
\left|\bar\partial\chi_z(\zeta)\right|_{i\partial\bar\partial\psi}^2
&=\left(\frac{4}{3\delta}\right)^2\left|\eta'\left(\frac{4}{3\delta} u_z(\zeta)\right)\right|^2
\left|\bar\partial u_z(\zeta)\right|_{i\partial\bar\partial\psi}^2\\
&\le C\frac{1}{\delta^2}\chi_z(\zeta).
\end{align*}
\end{proof}

There is a weighted $L^2$-estimate for the $\bar\partial$-problem due to Delin :
\begin{thm}[\cite{Del}, \cite{Hor}]\label{thm_Delin}
Assume that $\theta$ is a closed $(0,1)$-form on a pseudoconvex domain $\Omega \subset \Cn$ and 
that $\varphi$ is strictly plurisubharmonic and $\mathcal{C}^2$ there. 
Let $W$ be a weight function on $\Omega$ satisfying 
$|\partial W|_{i\partial\bar\partial\varphi}\leq \varepsilon W$ for some $\varepsilon \in (0,\sqrt{2}).$ 
Then the $L^2_\varphi(\Omega)$-minimal solution $f$ to the $\bar\partial$-equation
\begin{align*}
\overline{\partial}f=\theta
\end{align*}
satisfies
\begin{align*}
\int_{\Omega} |f|^2 e^{-\varphi} W\, dV
\leq \frac{2}{(\varepsilon - \sqrt{2})^2} \int_{\Omega} |\theta|^2_{i\partial\bar\partial\varphi} e^{-\varphi} W \,dV. 
\end{align*} 
Here, $|\theta|_{i\partial\bar\partial\varphi}$ denotes the norm of $\theta$ in the K\"{a}hler metric induced by the potential function $\varphi$. 
\end{thm} 

\subsection{Kernel estimates}
We have an estimate of the reproducing kernel on the diagonal. 

\begin{pro}\label{diagonal}
Let $z\in\Bn$. Then
\begin{align*}
K_\psi(z,z) \lesssim \frac{e^{\psi(z)}}{(1-|z|^2)^{2n+1}}.
\end{align*}
\end{pro}

\begin{proof}
For $z \in \Bn$, $K_z(w)$ is a holomorphic function on $\Bn$. 
By Lemma \ref{submean} for a small $\delta>0$, we get
\begin{align*}
|K_z(w)|^2 e^{-\psi(w)}
&\lesssim \frac{1}{(1-|w|^2)^{2n+1}} \int_{B_{\psi}(w,\delta)} |K_z(\zeta)|^2 e^{-\psi(\zeta)}dV(\zeta)\\
&\le \frac{1}{(1-|w|^2)^{2n+1}} \int_{\Bn} |K_z(\zeta)|^2 e^{-\psi(\zeta)}dV(\zeta)\\
&=\frac{K_z(z)}{(1-|w|^2)^{2n+1}}.
\end{align*} 
Hence we get 
\begin{align*}
|K_z(w)|^2\lesssim \frac{e^{\psi(w)}K_z(z)}{(1-|w|^2)^{2n+1}}.
\end{align*}  
By taking $w=z$, we obtain
\begin{align*}
K_{\psi}(z,z)\lesssim \frac{e^{\psi(z)}}{(1-|z|^2)^{2n+1}}.
\end{align*} 
\end{proof}

We prove our main result.
	\begin{thm}
	There are constants $C>0$ and $0<\varepsilon<\sqrt{2}$ satisfying
	\begin{align*}
	|K_{\psi}(z,w)|^2 \leq C \frac{e^{\psi(w)}e^{\psi(z)}}{(1-|w|^2)^{2n+1} (1-|z|^2)^{2n+1}} e^{-\varepsilon d_\psi(z, w)} 
	\quad \text{for} \quad z, w \in \Bn.
	\end{align*}
	\end{thm}

	\begin{proof}
	Let $\delta >0$ is fixed. First, we assume $d_\psi(z,w) \leq \delta$. 
	Then $1 \lesssim e^{-\varepsilon d_\psi(z,w)}$ for $\varepsilon>0$.
		Proposition \ref{diagonal} gives that
		\begin{align*}
		\left|K_\psi(z,w)\right|^2\leq K_\psi(z,z)K_\psi(w,w)
		\lesssim\frac{e^{\psi(z)+\psi(w)}}{(1-|z|^2)^{2n+1}(1-|w|^2)^{2n+1}}
		\end{align*}
		which implies

	\begin{align*}
	|K_{\psi}(z,w)|^2 \leq C \frac{e^{\psi(w)}e^{\psi(z)}}{(1-|w|^2)^{2n+1} (1-|z|^2)^{2n+1}} e^{-\varepsilon d_\psi(z, w)}.
	\end{align*}

	Next, we assume $d_\psi(z,w)>\delta$. 
	Let $\chi_z(\zeta)$ 
	be the characteristic function satisfying $\chi_z |_{B_{\psi}(z, \frac{\delta}{4})} \equiv 1$ and $\supp \chi_z \subset B_{\psi}(z, \delta)$
	in Lemma \ref{char_fn_2}.
	By Lemma \ref{submean}, it is obtained that
	\begin{align*}
	|K_z(w)|^2 e^{-\psi(w)} 
	&\lesssim \frac{1}{(1-|w|^2)^{2n+1}} \int_{B_\psi \left(w,\frac{\delta}{4}\right)} |K_z(\zeta)|^2 e^{-\psi(\zeta)}\, dV(\zeta)\\
	&\le \frac{1}{(1-|w|^2)^{2n+1}} \left\|K_z\right\|^2_{L^2(\chi_w e^{-\psi}dV)}.
	\end{align*}
	The norm of $K_z\in L^2(\chi_w e^{-\psi}dV)$ is given by
	\begin{align}\label{kernel_sup}
	\|K_z\|_{L^2(\chi_w e^{-\psi}dV)}=\sup_f\left| \langle f,K_z \rangle _{L^2(\chi_w e^{-\psi}dV)} \right|
	\end{align} 
	where $f$ is holomorphic on $B_\psi(w,\delta)$ with $\|f\|_{L^2(\chi_w e^{-\psi}dV)}=1$.
	Because $f\chi_w \in {L^2(\chi_w e^{-\psi}dV)}$, we have 
	\begin{align*}
	\langle f,K_z \rangle_{L^2(\chi_w e^{-\psi}dV)}=P_\psi(f\chi_w)(z)  
	\end{align*}
	where $P_\psi$ is the orthogonal projection to $A^2_\psi(\Bn)$. 

	Let $$v=(I-P_\psi)(\chi_w f),$$ then $v$ is the canonical solution of the equation $\overline{\partial}v=f \overline{\partial}\chi_w.$ Since $\chi_w(z)=0$, we have $|v(z)|=|P_\psi(\chi_w f)(z)|$. By Lemma \ref{submean},
	\begin{align*} 
	|P_\psi\chi_w f)(z)|^2 e^{-\psi(z)} 
	\lesssim\frac{1}{(1-|z|^2)^{2n+1}} \int_{B_{\psi}\left(z,\frac{\delta}{4}\right)} 
	\left|v\left(\zeta\right)\right|^2 e^{-\psi(\zeta)}\, dV(\zeta).
	\end{align*}
	We have
	\begin{align*} 
	|P_\psi(\chi_w f)(z)|^2 e^{-\psi(z)} 
	\lesssim \frac{1}{(1-|z|^2)^{2n+1}} \int_{B_{\psi}\left(z,\frac{\delta}{4}\right)} e^{-\varepsilon d_\psi(z, \zeta)} \left|v(\zeta)\right|^2 e^{-\psi(\zeta)}\, dV(\zeta)
	\end{align*}
	where the constant $\varepsilon>0$ is an undetermined positive number. 
	Let $u_z$ be the smoothing function of $d_\psi(z, \cdot)$ which is defined in the proof of Lemma \ref{char_fn_2}, then 
	\begin{align*} 
	|P_\psi(\chi_w f)(z)|^2 e^{-\psi(z)} 
	&\lesssim\frac{1}{(1-|z|^2)^{2n+1}} \int_{B_{\psi}\left(z,\frac{\delta}{4}\right)} e^{-\varepsilon u_z(\zeta)} 
	\left|v\left(\zeta\right)\right|^2 e^{-\psi(\zeta)}\, dV(\zeta)\\
	&\le\frac{1}{(1-|z|^2)^{2n+1}} \int_{\Bn} e^{-\varepsilon u_z(\zeta)} \left|v\left(\zeta\right)\right|^2 e^{-\psi(\zeta)} \,dV(\zeta). 
	\end{align*}

	Let $W(\zeta):=e^{-\varepsilon u_z(\zeta)}$, 
	then 
	$\left|\partial u_z\right|_{i\partial\bar\partial\psi}\le 1$ gives that
	\begin{align*}
	\partial W(\zeta) = -\varepsilon e^{-\varepsilon u_z(\zeta)} \partial u_z(\zeta)
	\end{align*}
	and
	\begin{align*}
	\left| \partial W\right|_{i\partial\bar\partial\psi}
	&=\varepsilon \left|\partial u_z\right|_{i\partial\bar\partial\psi} W\\
	&\leq \varepsilon W.
	\end{align*}
	We choose $\varepsilon>0$ sufficiently small so that $\varepsilon<\sqrt{2}$. 
	Since $\psi$ is strictly plurisubharmonic and $\left| \partial W\right|_{i\partial\bar\partial\psi} \leq \varepsilon W$, 
	Theorem \ref{thm_Delin} yields that 
	\begin{align*} 
	\int_{\Bn} \left|v\left(\zeta\right)\right|^2 e^{-\psi(\zeta)} W(\zeta)\,dV(\zeta) 
	&\lesssim\int_{\Bn} \left|f \overline{\partial}\chi_w\right|_{i\partial\bar\partial\psi}^2 e^{-\psi(\zeta)} W(\zeta)\,dV(\zeta)\\
	&=\int_{\Bn} \left|f(\zeta)\right|^2 \left|\overline{\partial}\chi_w(\zeta)\right|_{i\partial\bar\partial\psi}^2 e^{-\psi(\zeta)} W(\zeta)\,dV(\zeta).
	\end{align*}
	By Lemma \ref{char_fn_2}, we have $|\bar\partial {\chi_w}|_{i\partial\bar\partial\psi}^2 \leq C \frac{1}{\delta^2} \chi_w$ and so
	\begin{align*}
	\int_{\Bn} \left|v\left(\zeta\right)\right|^2 e^{-\psi(\zeta)} W(\zeta)\,dV(\zeta)
	&\lesssim\int_{\Bn} \left|f(\zeta)\right|^2 \chi_w(\zeta) e^{-\psi(\zeta)} W(\zeta)\,dV(\zeta)\\
	&=\int_{B_\psi(w,\delta)} \left|f(\zeta)\right|^2 e^{-\psi(\zeta)} W(\zeta)\,dV(\zeta).
	\end{align*}  
	Inequality (\ref{kerBn_smooth_fn_1}) yields that
	\begin{align*}
	\int_{\Bn} \left|v\left(\zeta\right)\right|^2 e^{-\psi(\zeta)} W(\zeta)\,dV(\zeta)
	\lesssim \int_{B_\psi(w,\delta)} \left|f(\zeta)\right|^2 e^{-\psi(\zeta)} e^{-\varepsilon d_\psi (z, \zeta)}\,dV(\zeta).
	\end{align*}  
	By triangle inequality, we have 
	$$
	d_\psi (z, \zeta) \leq d_\psi (z, w) + d_\psi (w, \zeta) \leq d_\psi (z, w) + \delta
	$$
	which implies that
	\begin{align*}
	\int_{\Bn} \left|v\left(\zeta\right)\right|^2 e^{-\psi(\zeta)} W(\zeta)\,dV(\zeta) 
	\lesssim e^{-\varepsilon d_\psi (z, w)} \int_{B_\psi(w,\delta)} \left|f(\zeta)\right|^2 e^{-\psi(\zeta)}\,dV(\zeta).
	\end{align*}  
	  
	For any holomorphic function $f$ on $B_\psi(w,\delta)$ with $\|f\|_{L^2(\chi_w e^{-\psi}dV)}=1$, it is obtained that
	\begin{align*}
	\int_{\Bn} \left|v\left(\zeta\right)\right|^2 e^{-\psi(\zeta)} W(\zeta)\,dV(\zeta) 
	\lesssim e^{-\varepsilon d_\psi(z, w)}.
	\end{align*} 
	Therefore, we have
	\begin{align*} 
	|P_\psi(\chi_w f)(z)|^2
	\lesssim \frac{e^{\psi(z)}}{(1-|z|^2)^{2n+1}} e^{-\varepsilon d_\psi(z, w)} 
	\end{align*}
	which gives the result
	\begin{align*}
	|K_z(w)|^2 \lesssim \frac{e^{\psi(w)}e^{\psi(z)}} {(1-|w|^2)^{2n+1} (1-|z|^2)^{2n+1}} e^{-\varepsilon d_\psi(z, w)}
	\end{align*}
	from (\ref{kernel_sup}). It completes the proof.
	\end{proof}

\section{Declarations}

\subsection{Ethical Approval}
Ethical approval not applicable to this article.

\subsection{Data Availability}
Data sharing not applicable to this article as no datasets were generated or analysed during
the current study.

\subsection{Funding}
H. R. Cho was supported by NRF  of Korea (NRF-2020R1F1A1A01048601)
and 
S. Park was supported by NRF of Korea (NRF-2021R1I1A1A01049889).

\subsection{Authors' contributions}
All the authors contributed to the writing of the present article. 
They also read and approved the final manuscript.

\subsection{Conflicts of interest}
On behalf of all authors, the corresponding author states that there is no conflict of interest.


\end{document}